\documentclass[10pt,oneside,english]{amsart}
\usepackage[T1]{fontenc}
\usepackage{geometry}
\usepackage{fancyhdr}
\usepackage{amssymb}

\makeatletter

\providecommand{\LyX}{L\kern-.1667em\lower.25em\hbox{Y}\kern-.125emX\@}

 \theoremstyle{plain}    
 \newtheorem{thm}{Theorem}[section]
 \numberwithin{equation}{section} 
 \numberwithin{figure}{section} 
 \usepackage{verbatim}
 \theoremstyle{plain}
 \theoremstyle{remark}    
 \newtheorem*{acknowledgement*}{Acknowledgement} 
 \theoremstyle{definition}
 \newtheorem{defn}[thm]{Definition}
 \theoremstyle{plain}    
 \newtheorem{prop}[thm]{Proposition} 
 \theoremstyle{definition}
  \newtheorem{example}[thm]{Example}
 \theoremstyle{plain}    
 \newtheorem{lem}[thm]{Lemma} 
 \theoremstyle{plain}    
 \newtheorem{cor}[thm]{Corollary} 

\usepackage{fancyhdr}

\def\quot{/\!\!/}
\def\tr{\mathsf{tr}}

\makeatother
\begin{document}

\title{Invariants of $2\times2$ matrices, irreducible $SL(2,\mathbb{C})$
characters and the Magnus trace map}

\author{Carlos A. A. Florentino}

\begin{abstract}
We obtain an explicit characterization of the stable points of the
action of $G=SL(2,\mathbb{C})$ on the cartesian product $G^{\times n}$
by simultaneous conjugation on each factor, in terms of the corresponding
invariant functions, and derive from it a simple criterion for irreducibility
of representations of finitely generated groups into $G$. We also
obtain analogous results for the action of $SL(2,\mathbb{C})$ on
the vector space of $n$-tuples of $2\times2$ complex matrices. For
a free group $F_{n}$ of rank $n$, we show how to generically reconstruct
the $2^{n-2}$ conjugacy classes of representations $F_{n}\to G$
from their values under the map $T_{n}:G^{\times n}\cong Hom(F_{n},G)\rightarrow\mathbb{C}^{3n-3}$
considered in \cite{M}, defined by certain $3n-3$ traces of words
of length one and two.
\end{abstract}
\maketitle

\section{Introduction and main results}

\baselineskip 4.4mmRepresentation varieties and character varieties
of finitely generated groups have been extensively studied in the
last three decades, not only for their many interesting properties,
but also in relation to subjects such as knot theory and spectral
geometry of hyperbolic manifolds, among several others (see for example
\cite{CS,Go,H,M} and references therein).

Here, we mainly concentrate on two problems related to the variety
of conjugacy classes of representations of finitely generated groups
into $G=SL(2,\mathbb{C})$, and particularly the case of representations
of a free group. The first is the characterization of the stable points
(in the sense of geometric invariant theory) of the action of $G$
on the cartesian product $G^{\times n}$ under simultaneous conjugation,
in terms of the corresponding invariant functions. As a consequence,
we obtain a simple numerical condition for the irreducibility of representations
of finitely generated groups. The second is a detailed study of a
map considered by Magnus \cite{M}, following earlier work by Vogt
\cite{V} and Fricke and Klein \cite{FK}, which is related to the
question of finding a minimal number of invariant functions on $G^{\times n}$
required to get all the other invariants by algebraic extensions.

We now describe the main results. Let $G$ be the algebraic Lie group
$SL(2,\mathbb{C})$ and, for a fixed integer $n\geq1$, let $X_{n}$
denote the cartesian product $G^{\times n}$. We are interested in
the orbit space of the action of $G$ on the affine variety $X_{n}$
under simultaneous conjugation on every factor. This is equivalent
to the space of conjugacy classes of $SL(2,\mathbb{C})$ representations
$\rho:F_{n}\rightarrow G$ of a free group $F_{n}$ on $n$ elements,
since\[
Hom(F_{n},G)\cong X_{n}\]
 by fixing a choice of generators of $F_{n}$.

In the context of algebraic geometry, we can consider the affine quotient
$X_{n}\quot G$, whose coordinate ring is the ring $\mathbb{C}[X_{n}]^{G}$
of regular functions on $X_{n}$ that are invariant under the action
of $G$. This is a categorical quotient where the geometric points
parametrize closed orbits. It is a consequence of a very general result,
the first fundamental theorem of invariants of $m\times m$ matrices
(see \cite{P} for $m\geq2$, or Thm. \ref{thm:Procesi} below for
$m=2$), that the building blocks of these $G$-invariant functions
are the following \emph{trace functions}. To any given word $w\in F_{n}$,
the corresponding trace function $t_{w}$ (sometimes called Fricke
character) is the algebraic $G$-invariant function\begin{equation}
t_{w}:Hom(F_{n},G)\rightarrow\mathbb{C}\label{eq:tracefct}\end{equation}
 that sends a representation $\rho$ to the trace of the $SL(2,\mathbb{C})$
matrix $\rho(w)$. Moreover, it is a very old result of Vogt and Fricke
(see \cite{V} and also \cite{H}) that the ring of these trace functions
is finitely generated. One of these finiteness results is as follows.
Let $\epsilon_{1},...,\epsilon_{n}\in F_{n}$ denote a fixed choice
of generators of $F_{n}$ and define the subset of $F_{n}$ consisting
of the lexicographically ordered words of length $\leq3$ with no
repeated letters \[
H_{n}=\{\epsilon_{j},\,1\leq j\leq n\}\cup\{\epsilon_{j}\epsilon_{k},\,1\leq j<k\leq n\}\cup\{\epsilon_{j}\epsilon_{k}\epsilon_{l},\,1\leq j<k<l\leq n\}\subset F_{n},\]
 of cardinality $N=n+\binom{n}{2}+\binom{n}{3}=\frac{n^{3}+5n}{6}$.
One can show that given any word $w\in F_{n}$, the function $t_{w}$
is a polynomial with rational coefficients in the variables $t_{\gamma}$,
$\gamma\in H_{n}$ (see \cite{V}, Cor. 4.14). These generators give
an embedding of the categorical quotient in $\mathbb{C}^{N}$, so
that $X_{n}\quot G$ corresponds to some polynomial ideal in $\mathbb{C}[t_{\gamma},\ \gamma\in H_{n}]$. 

In the present paper, inspired by geometric invariant theory (GIT),
we obtain a simple criterion, in terms of trace functions (\ref{eq:tracefct}),
for an element $A\in X_{n}$ to be in the subset $X_{n}^{\textrm{st}}$
of $X_{n}=G^{\times n}$ of stable points for the action of $G$.
By standard arguments of GIT, the \emph{affine stable quotient} $X_{n}^{\textrm{st}}/G$
will be an affine variety which is a geometric quotient of $X_{n}^{\textrm{st}}$
by $G$ in the sense that all fibers of the quotient map are indeed
orbits of the action. This is in contrast to the categorical affine
quotient $X_{n}\quot G$ where the points only parametrize closure-equivalence
classes of orbits. 

In another direction, we show that a similar numerical criterion can
be used to check the irreducibility of a representation of a finitely
generated group $\Gamma$ in $G=SL(2,\mathbb{C})$, as follows. Let
$\bar{\epsilon}_{1},...,\bar{\epsilon}_{n}$ be a choice of generators
of $\Gamma$. If $\rho:\Gamma\to G$ is a representation, define $A\in X_{n}$
by setting \[
A=(\rho(\bar{\epsilon}_{1}),...,\rho(\bar{\epsilon}_{n})).\]
 On the other hand, a given point $A=(A_{1},...,A_{n})\in X_{n}$
produces a representation $\rho$ of the free group $F_{n}$, by letting
$\rho(\epsilon_{j})=A_{j}$ for $j=1,...,n$. For three indices $1\leq j,k,l\leq n$,
denote by $\rho_{jkl}:F_{3}\to G$ the representation associated to
$(\rho(\bar{\epsilon}_{j}),\rho(\bar{\epsilon}_{k}),\rho(\bar{\epsilon}_{l}))\in X_{3}$.
In section 4, we show.

\begin{thm}
\label{thm:1}For a representation $\rho:\Gamma\rightarrow G$, the
following sentences are equivalent.

(i) $\rho:\Gamma\rightarrow G$ is reducible.

(ii) $(\rho(\bar{\epsilon}_{1}),...,\rho(\bar{\epsilon}_{n}))\in X_{n}$
is not stable for the conjugation action of $G$.

(iii) There is a $g\in SL(2,\mathbb{C})$ such that all matrices $g\rho(\bar{\epsilon}_{k})g^{-1}$
are upper triangular.

(iv) For all triples of indices $1\leq j,k,l\leq n$, $\rho_{jkl}:F_{3}\to G$
is reducible.
\end{thm}
In \cite{CS} (Cor. 1.2.2), Culler and Shalen proved that $\rho:\Gamma\to G$
is reducible if and only if $\mathsf{tr}(h)=2$ for every element
$h$ in the commutator subgroup $[\Gamma,\Gamma]$. Part (iv) of Theorem
\ref{thm:1} shows that irreducibility can be decided by a finite
process, looking only at all the associated representations of $F_{3}$.
Furthermore, the analysis of irreducibility for the case of $F_{3}$
leads to the following concrete numerical condition. For each triple
of indices $1\leq j,k,l\leq n$, define the following $G$ invariant
functions $\sigma_{jk},\Delta_{jkl}:X_{n}\to\mathbb{C}$ (see also
definition \ref{def:notation} below). \begin{eqnarray}
\sigma_{jk}(A) & = & \mathsf{tr}([A_{j},A_{k}])-2,\nonumber \\
\Delta_{jkl}(A) & = & (\mathsf{tr}(A_{j}A_{k}A_{l})-\mathsf{tr}(A_{l}A_{k}A_{j}))^{2},\label{eq:FD}\end{eqnarray}
where $[A_{j},A_{k}]=A_{j}A_{k}A_{j}^{-1}A_{k}^{-1}$ is the commutator
of $SL(2,\mathbb{C})$ matrices. The function $\Delta_{jkl}$ may
be called the \emph{Fricke discriminant}, being the discriminant of
the polynomial associated to the Fricke relation (see \cite{Go,M}).
We prove

\begin{thm}
\label{thm:2}Let $A=(A_{1},...,A_{n})\in X_{n}$ be the $n$-tuple
associated with the representation $\rho:\Gamma\to G$. Then $\rho$
is reducible if and only if $\sigma_{jk}(A)=\Delta_{jkl}(A)=0$ for
any triple $1\leq j,k,l\leq n$.
\end{thm}
The computations involved in the theorems above can be easily adapted
for the case of $SL(2,\mathbb{C})$ acting by simultaneous conjugation
on the vector space $V_{n}$ of $n$-tuples of arbitrary complex $2\times2$
matrices. In section 3, after briefly recalling the relevant definitions
in geometric invariant theory, we describe the stable locus for this
bigger space, and prove results analogous to the above theorem (see
Theorem \ref{thm:0}). All these results are based on the explicit
characterization of the $n$-tuples of $2\times2$ matrices that are
simultaneously similar to a set of $n$ upper triangular matrices,
in terms of invariant functions, which is obtained in section 2 (see
Theorem \ref{thm:UT}). In section 4, we also briefly comment on the
relation between this notion of stability and the stability of the
holomorphic vector bundle on a compact Riemann surface $S$ arising
from a representation of the fundamental group of $S$ into $G$ (see
Proposition \ref{pro:st-irr}).

Section 5 focus on the problem of reconstructing an orbit of the action
of $G$ on $X_{n}$ from a minimal number of traces, which was motivated
by the articles \cite{M} and \cite{Go}. Consider an arbitrary finite
sequence of words $J=(w_{1},...,w_{N})\in\left(F_{n}\right)^{N}$
and let $T_{J}$ denote the map\begin{eqnarray}
T_{J}:X_{n}\cong Hom(F_{n},G) & \rightarrow & \mathbb{C}^{N}\nonumber \\
\rho & \mapsto & (t_{w_{1}}(\rho),...,t_{w_{N}}(\rho)).\label{eq:Tj}\end{eqnarray}
Given that the quotient $G^{\times n}\quot G$ is a variety of dimension
$3n-3$, and $T_{J}$ factors through this quotient, it is natural
to look for a sequence $J$ of $N=3n-3$ words, such that $T_{J}$
is surjective onto a Zariski open subset of $\mathbb{C}^{3n-3}$ and
that all the preimages are finite (when non empty). Under the algebra-geometry
dictionary, this is equivalent to finding a minimal set of trace functions
$t_{w_{1}},...,t_{w_{N}}$ such that the field of invariant rational
functions on $G^{\times n}$ is an algebraic extension of $\mathbb{C}(t_{w_{1}},...,t_{w_{N}})$.

In this paper, we consider only those sequences $J$ composed of the
$n$ basic trace functions of length one $t_{\epsilon_{1}},...,t_{\epsilon_{n}}$
and some choice of $2n-3$ other words of length $2$. The basic example
is the result attributed to Vogt and Fricke that, for $n=2$, the
map\[
T_{(\epsilon_{1},\epsilon_{2},\epsilon_{1}\epsilon_{2})}:X_{2}\rightarrow\mathbb{C}^{3}\]
 is surjective. In \cite{Go}, Goldman presents an almost self-contained
proof of this, showing also that for $n=3$, the trace map \[
T_{(\epsilon_{1},\epsilon_{2},\epsilon_{3},\epsilon_{1}\epsilon_{2},\epsilon_{1}\epsilon_{3},\epsilon_{2}\epsilon_{3})}:X_{3}\rightarrow\mathbb{C}^{6}\]
 is again surjective. To consider the case of a free group of arbitrary
rank $n\geq4$, first note that we cannot take as $J$ the sequence
with the $n$ words of length one together with the $\binom{n}{2}$
ordered words of length two, since this would have length greater
than the wanted $3n-3$. Let us choose the sequence of length $3n-3$
omitting those words of length two $\epsilon_{j}\epsilon_{k}$ such
that $\{ j,k\}\cap\{1,2\}=\emptyset$, \begin{equation}
J_{n}:=(\epsilon_{1},\epsilon_{2},\epsilon_{1}\epsilon_{2},...,\epsilon_{k},\epsilon_{1}\epsilon_{k},\epsilon_{2}\epsilon_{k},...,\epsilon_{n},\epsilon_{1}\epsilon_{n},\epsilon_{2}\epsilon_{n}),\label{eq:J}\end{equation}
and denote the corresponding trace map by $T_{n}=T_{J_{n}}$. Of course,
it would be completely analogous to use another sequence of the form
(\ref{eq:J}) with another pair of indices playing the role of $\{1,2\}$.
In terms of the $n$-tuple $A=(A_{1},...,A_{n})\in G^{\times n}$
corresponding to $\rho:F_{n}\to G$ we have\begin{equation}
T_{n}(A)=(t_{1},t_{2},t_{12},...,t_{k},t_{1k},t_{2k},...,t_{n},t_{1n},t_{2n}),\label{eq:T}\end{equation}
where we use the notation $t_{j}=\mathsf{tr}(A_{j})$, $t_{jk}=\mathsf{tr}(A_{j}A_{k})$.
In section 5, we show that $T_{n}$ is almost surjective for $n\geq4$,
and omits a set contained in a very explicit irreducible subvariety
of $\mathbb{C}^{3n-3}$ of codimension $1$. Examples of points $\mathbf{z}\in\mathbb{C}^{3n-3}$,
with $T_{n}^{-1}(\mathbf{z})$ empty are given in the Appendix. 

The algebraic map $T_{n}:X_{n}\to\mathbb{C}^{3n-3}$ (\ref{eq:T})
will be called the \emph{Magnus trace map}. In \cite{M}, Magnus showed
that, given three matrices $A_{1},A_{2}$ and $A_{3}$ verifying $\sigma_{12}(A_{1},A_{2})\neq0$
and $\Delta_{123}(A_{1},A_{2},A_{3})\neq0$ and given any $\mathbf{q}\in\mathbb{C}^{3n-9}$
(thought as the last $3n-9$ coordinates in $\mathbb{C}^{3n-3}$)
one can find $n-3$ other matrices $A_{4},...,A_{n}$ such that $T_{n}(A)=(\mathbf{p},\mathbf{q})$,
where $\mathbf{p}=(t_{1},t_{2},t_{12},t_{3},t_{13},t_{23})$. He also
proved that the number of different solutions $(A_{4},...,A_{n})$
realizing this equation is bounded by $2^{n-2}$. It turns out that
the Fricke discriminant condition $\Delta_{123}\neq0$ is not really
necessary, and we only need to impose the condition $\sigma_{12}(A_{1},A_{2})\neq0$
to be able to find $A_{3},...,A_{n}$ such that $\mathsf{tr}(A_{k})$,
$\mathsf{tr}(A_{1}A_{k})$ and $\mathsf{tr}(A_{2}A_{k})$ assume preassigned
values for $k=3,...,n$. Moreover, the condition $\sigma_{12}(A_{1},A_{2})\neq0$
allows us to express the (at most $2^{n-2}$) orbits in the preimage
of $T_{n}$ in very explicit terms.

\begin{thm}
\label{thm:Magnus}Let $B_{1},B_{2}\in SL(2,\mathbb{C})$ be such
that $\mathsf{tr}([B_{1},B_{2}])\neq2$. Then, there exists a $g\in G$
such that $A_{j}:=gB_{j}g^{-1}$, $j=1,2$ are invariant under transposition.
Let $n\geq3$ and $\mathbf{r}=(t_{1},t_{2},t_{12})$. Then, given
any $\mathbf{s}\in\mathbb{C}^{3n-6}$ there exist $n-2$ matrices
$A_{3},...,A_{n}\in SL(2,\mathbb{C})$ such that \begin{equation}
T_{n}(A)=(\mathbf{r},\mathbf{s}).\label{eq:trace}\end{equation}
Moreover, given any solution $A\in X_{n}$ of (\ref{eq:trace}) with
$\mathsf{tr}([A_{1},A_{2}])\neq2$ and $A_{1}^{T}=A_{1}$ and $A_{2}^{T}=A_{2}$
(where $T$ denotes transposition) the inverse image $T_{n}^{-1}(\mathbf{r},\mathbf{s})$
consists of the $G$ orbits of the finite set\[
\left\{ (A_{1},A_{2},B_{3},...,B_{n}):B_{j}=A_{j}\textrm{ or }B_{j}=A_{j}^{T}\textrm{ for }j=3,...,n\right\} .\]

\end{thm}
We would like to mention that most of the methods in this article
are constructive, in the sense that they can be used to implement
algorithms to effectively compute the quantities involved.

\begin{acknowledgement*}
We thank W. Goldman for his interest and encouragement, and my colleagues
J. Mourão and J. P. Nunes for many interesting and motivating conversations
on this and related topics. This work was partially supported by Center
for Mathematical Analysis, Geometry and Dynamical Systems, IST, and
by the {}``Fundação para a Ciência e a Tecnologia'' through the
programs Praxis XXI, POCI/MAT/58549/2004 and FEDER. Typeset using
\LyX{}. 
\end{acknowledgement*}

\section{Degenerate simultaneous similarity of $2\times2$ matrices}

In this section, we are interested in the simultaneous conjugacy classes
of a finite set of $2\times2$ complex matrices. We will describe
the most degenerate cases, in particular give necessary and sufficient
conditions, in terms of invariant functions, for $n$ matrices to
be simultaneously conjugated to matrices in upper or lower triangular
form.

Let the general linear group $GL(2,\mathbb{C})$ act on the vector
space of $n$-tuples of $2\times2$ matrices ($n\geq1$)\[
V_{n}:=\left(M_{2\times2}(\mathbb{C})\right)^{\times n}\]
 by simultaneous conjugation\begin{equation}
g\cdot A:=(gA_{1}g^{-1},...,gA_{n}g^{-1}),\label{accao}\end{equation}
 where $A=(A_{1},...,A_{n})\in V_{n}$ and $g\in GL(2,\mathbb{C})$.
There are plenty of $GL(2,\mathbb{C})$-invariant regular (i.e, polynomial)
functions on $V_{n}$ and, by the first fundamental theorem of invariants
of $m\times m$ matrices \cite{P}, the trace functions defined by\[
V_{n}\rightarrow\mathbb{C},\quad\quad A\mapsto\tr(A_{i_{1}}\cdots A_{i_{k}})\]
 and labelled by `words' $A_{i_{1}}\cdots A_{i_{k}}$ in the components
$A_{j}$ of $A\in V_{n}$, generate the ring of invariants $\mathbb{C}[V_{n}]^{GL(2,\mathbb{C})}$.
Moreover, this ring is finitely generated and we have

\begin{thm}
\emph{(Procesi \cite{P})} \label{thm:Procesi}Any $GL(2,\mathbb{C})$-invariant
regular function on $V_{n}$ is a polynomial in the following set
of $(n+1)^{3}$ functions\[
A\mapsto\tr(A_{j}A_{k}A_{l}),\quad\quad0\leq j,k,l\leq n,\]
 where $A=(A_{1},...,A_{n})\in V_{n}$ and $A_{0}=I$ is the identity
$2\times2$ matrix.
\end{thm}
The many relations between these functions, described by the second
fundamental theorem of invariants of matrices (see \cite{P}), will
not be important here. 

As in the case $n=1$, two elements $A$ and $A'$ of $V_{n}$ will
be called \emph{similar} if they are in the same $GL(2,\mathbb{C})$
orbit. Note that an element $A\in V_{n}$ can be viewed either as
a vector of $2\times2$ matrices as above or, alternatively, as a
\emph{single} matrix with vector valued entries\begin{equation}
A=\left(\begin{array}{cc}
a & b\\
c & d\end{array}\right),\quad a,b,c,d\in\mathbb{C}^{n}.\label{eq:nmatrix}\end{equation}
 This justifies the following terminology and notation.

\begin{defn}
A point $A=(A_{1},...,A_{n})\in V_{n}$ will be called an $n$-matrix.
We will say that $A$ is an \emph{upper triangular $n$-matrix} if
the vector $c\in\mathbb{C}^{n}$ is zero, and we will denote by $\mathcal{UT}_{n}\subset V_{n}$
the $GL(2,\mathbb{C})$ orbit of the subset of upper triangular $n$-matrices.
Hence, $A\in\mathcal{UT}_{n}$ if and only if $A$ is \emph{similar}
to an upper triangular $n$-matrix. $A_{j}$ is called the $j$th
\emph{component} of $A$. 
\end{defn}
Note that $A$ is similar to an upper triangular $n$-matrix if and
only if there is a proper nonzero subspace of $\mathbb{C}^{2}$ which
is fixed by every component $A_{j}$ of $A$. The condition $A\in\mathcal{UT}_{n}$
is also equivalent to $A$ being similar to a lower triangular $n$-matrix
(one with zero $b\in\mathbb{C}^{n}$).

The similarity classes of a pair of $m\times m$ matrices were obtained
in \cite{Fr}, and in the simplest $m=2$, $n=2$ case, the following
irreducible algebraic subset of $V_{2}$ plays an important role\begin{equation}
W=\left\{ (A_{1},A_{2})\in V_{2}:(t_{11}-\frac{1}{2}t_{1}^{2})(t_{22}-\frac{1}{2}t_{2}^{2})=\left(t_{12}-\frac{1}{2}t_{1}t_{2}\right)^{2}\right\} .\label{s12}\end{equation}
 Here and below, we are using the following notation \begin{eqnarray*}
t_{j} & := & \tr(A_{j})\\
t_{jk} & := & \tr(A_{j}A_{k})\end{eqnarray*}
 for a general element $A=(A_{1},...,A_{n})\in V_{n}$ and any pair
of indices $j,k\in\{1,...,n\}$. The relevance of $W$ can be seen
from the fact that if $(A_{1},A_{2})$ does not belong to $W$, then
its $GL(2,\mathbb{C})$ orbit is uniquely determined by $t_{1},t_{2},t_{11},t_{22},t_{12}$
(\emph{\cite{Fr}}). Moreover, the following is not difficult to prove.

\begin{prop}
\emph{(see, for instance, \cite{Fr})}\label{pair} A pair $(A_{1},A_{2})\in V_{2}$
belongs to $W$ if and only if it is in the orbit of a pair of upper
triangular matrices.
\end{prop}
To generalize this result to higher $n$, let us abbreviate some frequently
used $GL(2,\mathbb{C})$-invariant functions on $V_{n}$ as follows

\begin{defn}
\label{def:notation}Define, for every triple of indices $1\leq j,k,l\leq n$,\begin{eqnarray*}
\tau_{jk} & := & t_{jk}-\frac{1}{2}t_{j}t_{k},\\
\sigma_{jk} & := & \tau_{jk}^{2}-\tau_{jj}\tau_{kk},\\
\Delta_{jkl} & := & (t_{jkl}-t_{lkj})^{2}.\end{eqnarray*}
 We omit the dependence on $A=(A_{1},...,A_{n})$ where no ambiguity
arises. For $A\in V_{n}$ and $g\in GL(2,\mathbb{C})$, we will always
write\[
A_{j}=\left(\begin{array}{cc}
a_{j} & b_{j}\\
c_{j} & d_{j}\end{array}\right),\qquad j=1,...,n,\qquad g=\left(\begin{array}{cc}
x & y\\
z & w\end{array}\right)\]
or sometimes $A_{j}=(a_{j},b_{j},c_{j},d_{j})$, $g=(x,y,z,w)$. In
terms of these variables and of $e_{j}:=a_{j}-d_{j}$, we have:\begin{eqnarray}
\tau_{jk} & = & \frac{e_{j}e_{k}}{2}+b_{j}c_{k}+c_{j}b_{k},\nonumber \\
\sigma_{jk} & = & (b_{j}c_{k}-c_{j}b_{k})^{2}-(b_{j}e_{k}-e_{j}b_{k})(c_{j}e_{k}-e_{j}c_{k}).\label{eq:s}\end{eqnarray}
 We also use the abbreviation\[
\nu_{j}:=\tau_{jj}=t_{jj}-\frac{1}{2}t_{j}^{2}=\frac{e_{j}^{2}}{2}+2b_{j}c_{j},\]
so that $\nu_{j}\neq0$ if and only if $A_{j}$ has distinct eigenvalues.
The functions $\sigma_{jk}$ and $\Delta_{jkl}$ are fully symmetric
under index permutation and vanish whenever two indices coincide.
Note that $-\sigma_{jk}$ and $-\frac{1}{2}\Delta_{jkl}$ are, respectively,
the top left $2\times2$ minor and the determinant of the symmetric
matrix\begin{equation}
\left(\begin{array}{ccc}
\tau_{jj} & \tau_{jk} & \tau_{jl}\\
\tau_{kj} & \tau_{kk} & \tau_{kl}\\
\tau_{lj} & \tau_{lk} & \tau_{ll}\end{array}\right).\label{eq:det}\end{equation}

\end{defn}
With slightly different normalizations, the restrictions of these
functions to $SL(2,\mathbb{C})^{\times n}$ were used in \cite{GM}
and \cite{M}. Since the equation (\ref{s12}) that defines $W\subset V_{2}$
is $\sigma_{12}=0$, the condition $\sigma_{jk}(A)=0$ for some $j,k\in\{1,...,n\}$
is equivalent to $(A_{j},A_{k})\in\mathcal{UT}_{2}$, by Proposition
\ref{pair}. Therefore, we have

\begin{prop}
\label{UT->s0}If $A\in V_{n}$ is similar to an upper triangular
$n$-matrix, then $\sigma_{jk}(A)=0$ for all distinct $1\leq j,k\leq n$. 
\end{prop}
\begin{proof}
If $A=(A_{1},...,A_{n})\in\mathcal{UT}_{n}$ then $g\cdot A$ is upper
triangular, for some $g\in G$. Hence, $g\cdot(A_{j},A_{k})$ is an
upper triangular $2$-matrix, for any $1\leq j,k\leq n$, and $\sigma_{jk}=0$
by Proposition \ref{pair}.
\end{proof}
For $n\geq2$, define, for distinct $j,k\in\{1,...,n\}$, the following
algebraic subsets of $V_{n}$\[
W_{jk}=W_{kj}:=\{ A\in V_{n}:\sigma_{jk}(A)=0\},\]
 \[
\Sigma_{n}:=\bigcap_{j,k}W_{jk}.\]
Then, $A\in\Sigma_{n}$ if and only if every pair of matrix components
of $A$ is in $\mathcal{UT}_{2}$. From Proposition \ref{UT->s0},
we have\[
\mathcal{UT}_{n}\subset\Sigma_{n}\subset V_{n}\]
 for all $n\geq2$. However, the vanishing of all $\sigma_{jk}$ is
not sufficient for $A$ to be in $\mathcal{UT}_{n}$, for $n\geq3$,
as the next example shows. 

\begin{example}
\label{exa:example}Let $A=(A_{1},A_{2},A_{3})$ be defined by\[
A_{1}=\left(\begin{array}{cc}
a_{1} & 0\\
0 & d_{1}\end{array}\right),\quad A_{2}=\left(\begin{array}{cc}
a_{2} & b_{2}\\
0 & d_{2}\end{array}\right),\quad A_{3}=\left(\begin{array}{cc}
a_{3} & 0\\
c_{3} & d_{3}\end{array}\right).\]
 Then $\sigma_{12}=\sigma_{13}=0$ and $\sigma_{23}=b_{2}c_{3}(e_{2}e_{3}+b_{2}c_{3})$.
Assume that $e_{2}e_{3}+b_{2}c_{3}=0$ and that $e_{1}b_{2}c_{3}\neq0$,
so that all $\sigma_{jk}$ vanish, neither $A_{2}$ or $A_{3}$ are
diagonal, and (since these assumptions imply $e_{2}e_{3}\neq0$) all
three matrices have distinct eigenvalues. Let $g=(x,y,z,w)\in SL(2,\mathbb{C})$.
Then\begin{eqnarray}
gA_{1}g^{-1} & = & \left(\begin{array}{cc}
* & -xye_{1}\\
zwe_{1} & *\end{array}\right)\nonumber \\
gA_{2}g^{-1} & = & \left(\begin{array}{cc}
* & x(xb_{2}-ye_{2})\\
z(we_{2}-zb_{2}) & *\end{array}\right)\label{eq:BT}\\
gA_{3}g^{-1} & = & \left(\begin{array}{cc}
* & -y(xe_{3}+yc_{3})\\
w(ze_{3}+wc_{3}) & *\end{array}\right),\nonumber \end{eqnarray}
from which it follows that $A$ is not similar to an upper triangular
$3$-matrix. 
\end{example}
On the other hand, we have.

\begin{thm}
\label{thm:UT} Let $n\geq1$. $A=(A_{1},...,A_{n})\in V_{n}$ is
similar to an upper triangular matrix if and only if for all triples
$1\leq j,k,l\leq n$, $(A_{j},A_{k},A_{l})$ is similar to an upper
triangular $3$-matrix.
\end{thm}
For the proof, we will use the following lemmata.

\begin{lem}
\label{lemma}Let $A$ be an upper triangular $n$-matrix with $\nu_{j}=\frac{e_{j}^{2}}{2}\neq0$,
for some $1\leq j\leq n$. Then, $A$ is similar to another upper
triangular $n$-matrix $A'=(A_{1}',...,A_{n}')$ with $A_{j}'$ diagonal.
\end{lem}
\begin{proof}
We only need to find $g\in GL(2,\mathbb{C})$ such that $gA_{k}g^{-1}$
is still upper triangular for any $k$, and such that $gA_{j}g^{-1}$
is diagonal. Letting $g=(x,y,0,x^{-1})$ for some $x\neq0$, we calculate
$gA_{k}g^{-1}=(a_{k},b_{k}x^{2}-yxe_{k},0,d_{k})$, for every $k=1,...,n$.
Therefore, all $gA_{k}g^{-1}$ are upper triangular matrices, and
using $y=b_{j}x/e_{j}$ (since $e_{j}\neq0$), $gA_{j}g^{-1}$ is
diagonal. 
\end{proof}
\begin{lem}
\label{lem:lemma2}As in Example \ref{exa:example}, let $A=(A_{1},A_{2},A_{3})$
be a triple of the form\begin{equation}
A_{1}=\left(\begin{array}{cc}
a_{1} & 0\\
0 & d_{1}\end{array}\right),\quad A_{2}=\left(\begin{array}{cc}
a_{2} & b_{2}\\
0 & d_{2}\end{array}\right),\quad A_{3}=\left(\begin{array}{cc}
a_{3} & 0\\
c_{3} & d_{3}\end{array}\right),\label{eq:form}\end{equation}
 Then $A\in\mathcal{UT}_{3}$ if and only if $e_{1}b_{2}c_{3}=0$. 
\end{lem}
\begin{proof}
If $e_{1}b_{2}c_{3}=0$ at least one of the factors is zero. In each
case, $A_{1}$ is a scalar, $A$ is lower triangular, or $A$ is upper
triangular, respectively, so $A\in\mathcal{UT}_{3}$. Conversely,
suppose that $e_{1}b_{2}c_{3}\neq0$. Then, from equations (\ref{eq:BT})
we see that there is no $g\in G$ that will make $g\cdot A$ upper
or lower triangular, so $A\notin\mathcal{UT}_{3}$.
\end{proof}
We can now finish the proof of Theorem \ref{thm:UT}.

\begin{proof}
The Theorem is obvious for $n\leq3$, so let $n\geq4$. If $A\in V_{n}$
is similar to an upper triangular $n$-matrix, then obviously any
$m$-tuple $(A_{j_{1}},...,A_{j_{m}})$ of $m\leq n$ components of
$A$ will be in $\mathcal{UT}_{m}$. Conversely, let all triples of
components be in $\mathcal{UT}_{3}$ and suppose, by induction, that
the result is valid for $n-1$. Then, in particular, all the $\sigma_{jk}$
and all $\Delta_{jkl}$ are zero, for indices $j,k,l$ between $1$
and $n-1$. To reach a contradiction, assume that $A$ is not similar
to an upper triangular $n$-matrix. By induction, we can suppose that
$(A_{1},...,A_{n-1})$ has been conjugated so that it is already an
upper triangular $(n-1)$-matrix. Let $A_{n}=(a_{n},b_{n},c_{n},d_{n})$
with $c_{n}\neq0$. None of the $A_{j}$ can be central, otherwise
the result would follow by induction. The $n-1$ conditions $\sigma_{jn}=0$,
$j=1,...,n-1$ imply (because $c_{n}\neq0$)\begin{eqnarray}
b_{j}^{2}c_{n}+b_{j}e_{j}e_{n}-e_{j}^{2}b_{n} & = & 0,\quad\textrm{for }j=1,...,n-1.\label{eq:inds0}\end{eqnarray}
 If one of the $e_{j}$, $j=1,...,n-1$ vanishes then, since $A_{j}$
is non-scalar, $b_{j}\neq0$ and the equations (\ref{eq:inds0}) become
$b_{j}^{2}c_{n}=0$ and have no solution. As a consequence, none of
these $e_{j}$'s can be zero. Then, by lemma \ref{lemma}, one can
assume that $b_{1}=0$, and the equation (\ref{eq:inds0}) with $j=1$
implies $b_{n}=0$. Now, we have all triples $(A_{1},A_{k},A_{n})$,
with $k=2,...,n-1$, in the form (\ref{eq:form}). Since $(A_{1},A_{k},A_{n})\in\mathcal{UT}_{3}$
by hypothesis, lemma \ref{lem:lemma2} implies $e_{1}b_{k}c_{n}=0$.
So, $b_{k}=0$ for all $k=2,...,n-1$. But then $(A_{1},...,A_{n})$
is lower triangular, and we have a contradiction.
\end{proof}
Note that the Theorem is true for any algebraically closed field of
characteristic 0 in place of the field of complex numbers. The following
statement is also useful.

\begin{lem}
\label{lem:S-UT}Let $A=(A_{1},A_{2},A_{3})\in\Sigma_{3}\setminus\mathcal{UT}_{3}$.
Then, $A$ is similar to a triple of the form (\ref{eq:form}) with
$e_{2}e_{3}+b_{2}c_{3}=0$ and $e_{1}b_{2}c_{3}\neq0$.
\end{lem}
\begin{proof}
Assume that $A\in\Sigma_{3}$. In particular, $\sigma_{12}=0$, so
we can suppose that $A_{1}$ and $A_{2}$ are both upper triangular.
Let $A_{3}=(a_{3},b_{3},c_{3},d_{3})$. Since $A\notin\mathcal{UT}_{3}$,
none of the $A_{j}$ can be a scalar, and $c_{3}$ is nonzero. The
2 conditions $\sigma_{j3}=0$, $j=1,2$ imply (because $c_{3}\neq0$)\begin{eqnarray}
b_{j}^{2}c_{3}+b_{j}e_{j}e_{3}-e_{j}^{2}b_{3} & = & 0,\quad\textrm{for }j=1,2.\label{eq:s0}\end{eqnarray}
 If one of the $e_{j}$, $j=1,2$ vanishes then, since $A_{j}$ is
non-scalar, $b_{j}\neq0$ and the equations (\ref{eq:s0}) become
$b_{j}^{2}c_{3}=0$ and have no solution. So, necessarily $e_{1}e_{2}\neq0$.
Then, by lemma \ref{lemma}, one can assume that $b_{1}=0$, and the
equation (\ref{eq:s0}) with $j=1$ implies $b_{3}=0$. Since $A$
cannot be lower triangular, we need to have $b_{2}\neq0$, and (\ref{eq:s0})
for $j=2$ simplifies to $b_{2}c_{3}+e_{2}e_{3}=0$.
\end{proof}
\begin{cor}
$(A_{1},A_{2},A_{3})\in V_{3}$ is similar to an upper triangular
$3$-matrix if and only if $\sigma_{12}=\sigma_{13}=\sigma_{23}=\Delta_{123}=0$.
\end{cor}
\begin{proof}
If $A\in\mathcal{UT}_{3}$, then $\sigma_{jk}=0$ for all $j,k$,
by Proposition \ref{UT->s0} and by direct computation $t_{123}=t_{321}$.
Conversely, if all $\sigma_{jk}=0$ and $A\notin\mathcal{UT}_{3}$,
by lemma \ref{lem:S-UT}, we can suppose that $A$ is in the form
(\ref{eq:form}) with $e_{1}b_{2}c_{3}\neq0$. An easy calculation
then gives $\Delta_{123}=(t_{123}-t_{321})^{2}=e_{1}^{2}b_{2}^{2}c_{3}^{2}\neq0$.
\end{proof}
This Corollary extends to all $2\times2$ matrices the result stated
in Prop. 4.4 of \cite{GM} for triples of $SL(2,\mathbb{C})$ matrices.
The following Proposition will be useful later.

\begin{prop}
\label{pro:ST}Let $A$ be a $2$-matrix with $\sigma_{12}\neq0$.
Then $A$ is similar to a $2$-matrix $B$ invariant under simultaneous
transposition ($B_{1}^{T}=B_{1}$ and $B_{2}^{T}=B_{2}$).
\end{prop}
\begin{proof}
We consider first the case when at least one of the matrices $A_{1}$
or $A_{2}$ is diagonalizable. Without loss of generality let $A_{1}$
be diagonal and $A_{2}=\left(a_{2},b_{2},c_{2},d_{2}\right)$. Then
$\sigma_{12}=-e_{1}^{2}b_{2}c_{2}$ from (\ref{eq:s}). Conjugation
of $A$ by the diagonal matrix $g=(x,0,0,x^{-1})$ produces the assignment
$b_{2}\mapsto b_{2}x^{2}$, $c_{2}\mapsto c_{2}x^{-2}$ and so, we
just solve $b_{2}x^{2}=c_{2}x^{-2}$ for an appropriate $x\neq0$,
which is possible since $\sigma_{12}\neq0$. Now, consider the case
$\nu_{1}=\nu_{2}=0$. Assuming that $A_{1}$ is already in Jordan
canonical form, write $A=(A_{1},A_{2})$ as \begin{equation}
A_{1}=\left(\begin{array}{cc}
a_{1} & b_{1}\\
0 & a_{1}\end{array}\right),\qquad A_{2}=\left(\begin{array}{cc}
a_{2} & b_{2}\\
c_{2} & d_{2}\end{array}\right),\label{eq:vpi}\end{equation}
with $b_{1}c_{2}\neq0$ according to the hypothesis $\sigma_{12}\neq0$.
Suppose first that $b_{2}=0$, which implies $d_{2}=a_{2}$, since
$\nu_{2}=0$. Conjugating $(A_{1},A_{2})$ by a diagonal matrix as
before, we can further assume that $b_{1}=c_{2}$. Then, using \[
g=\frac{1}{\sqrt{2}}\left(\begin{array}{cc}
1 & i\\
i & 1\end{array}\right),\]
and we obtain explicitly the transposition invariant pair,\[
B_{1}=g\cdot A_{1}=\left(\begin{array}{cc}
a_{1}+\lambda & i\lambda\\
i\lambda & a_{1}-\lambda\end{array}\right),\qquad B_{2}=g\cdot A_{2}=\left(\begin{array}{cc}
a_{2}-\lambda & i\lambda\\
i\lambda & a_{2}+\lambda\end{array}\right),\]
with $b_{1}=c_{2}=2i\lambda$. Finally, if $b_{2}\neq0$ in (\ref{eq:vpi}),
using the equation $\nu_{2}=\frac{e_{2}^{2}}{2}+2b_{2}c_{2}=0$, which
implies $e_{2}\neq0$, it is a simple computation to show that conjugation
by $g=(x,1,0,x^{-1})$ with $x=\frac{e_{2}}{2b_{2}}$ reduces that
pair to one with $A_{1}$ upper and $A_{2}$ lower triangular. Note
that in both cases, the transposition invariant $2$-matrix $B$ verifies
$\sigma_{12}=-e_{1}^{2}b_{2}^{2}=16\lambda^{4}$.
\end{proof}

\section{Stability of the $SL(2,\mathbb{C})$ action on $2\times2$ matrices}

In this section, we determine, in terms of invariant functions, the
stable points of the action of $SL(2,\mathbb{C})$ on the vector space
$V_{n}$ of $n$-tuples of complex $2\times2$-matrices, under simultaneous
conjugation on each factor.

Recall that, in the general situation of a general algebraic linearly
reductive Lie group $K$ acting on an affine variety $V$ one defines
the affine quotient variety $V\quot K$ as the spectrum of the ring
of invariant functions on $V$, which comes equipped with a projection\[
q:V\rightarrow V\quot K\]
 induced from the canonical inclusion of algebras $\mathbb{C}[V]\subset\mathbb{C}[V]^{K}$.
The set of closed orbits is in bijective correspondence with geometric
points of the quotient $V\quot K$. Recall also that a vector $x\in V$
is said to be \emph{stable} if the corresponding map\[
\psi_{x}:K\rightarrow V,\quad\quad g\mapsto g\cdot x\]
 is proper. It is easy to see that $x\in V$ is stable if and only
if the closure of the $K$-orbit of $x$ does not intersect the closed
subset consisting of points $x\in V$ with positive dimensional stabilizer
subgroup. Another useful criterion for stability is the \emph{Hilbert-Mumford
numerical criterion}, which is stated in terms of nontrivial homomorphisms
$\phi:\mathbb{C}^{*}\rightarrow K$, called one parameter subgroups
(1PS) of $K$. To any such $\phi$ and to a point $x\in V$ one associates
the morphism $\phi_{x}:\mathbb{C}^{*}\rightarrow V$ given by mapping
$\lambda\in\mathbb{C}^{*}$ to the point $\phi(\lambda)\cdot x$.
If $\phi_{x}$ can be extended to a morphism $\overline{\phi_{x}}:\mathbb{C}\rightarrow V$,
we say that $\lim_{\lambda\to0}\phi_{x}$ exists and equals $\overline{\phi_{x}}(0)$. 

\begin{thm}
\emph{(Hilbert-Mumford \cite{MFK}, see also \cite{Gi})} \label{thm:HMC}A
point $x\in V$ is stable if and only for every one parameter subgroup
$\phi$ of $K$, $\phi_{x}$ cannot be extended to a morphism $\mathbb{C}\rightarrow V$.
\end{thm}
It is easy to see that the conjugation action of $GL(2,\mathbb{C})$
on the space of $n$-tuples of $2\times2$ matrices $V_{n}$ has no
stable points, since the scalar nonzero matrices will stabilize any
point $A\in V_{n}$. This is not a big problem, since the same orbit
space can be obtained with the conjugation action of $G\equiv SL(2,\mathbb{C})$
on $V_{n}$ which has generically finite stabilizers. This is just
the restriction to $G\subset GL(2,\mathbb{C})$ of the action (\ref{accao}).
One could as well consider the action of $PSL(2,\mathbb{C})$ which
would have generically trivial stabilizers, but we will keep using
$G=SL(2,\mathbb{C})$. It is clear that any diagonal $n$-matrix $A$
(one for which both vectors $b$ and $c$ in (\ref{eq:nmatrix}) are
zero) has the subgroup \begin{equation}
H=\left(\begin{array}{cc}
\lambda & 0\\
0 & \lambda^{-1}\end{array}\right)\subset G,\,\,\lambda\in\mathbb{C}^{*}\label{eq:subg}\end{equation}
 contained in its stabilizer. The following is then an easy application
of the Hilbert-Mumford criterion, Theorem \ref{thm:HMC}.

\begin{prop}
\label{pro:stable}$A\in V_{n}$ is stable if and only if $A$ is
not similar to an upper triangular $n$-matrix.
\end{prop}
\begin{proof}
If $A$ is an upper triangular $n$-matrix, a simple computation shows
that the closure of the orbit of $A$ under the subgroup $H\subset G$
(\ref{eq:subg}) will intersect $D$. Therefore, no point in the orbit
of $A$ will be stable. Conversely, let $A\in V_{n}$ be not stable
and apply the numerical criterion. By elementary representation theory,
any one parameter subgroup of $G$ is conjugated to\[
\lambda\mapsto\phi_{n}(\lambda)=\left(\begin{array}{cc}
\lambda^{n} & 0\\
0 & \lambda^{-n}\end{array}\right),\quad n\in\mathbb{N}_{0}.\]
 In other words, any 1PS can be written as $\phi=g^{-1}\phi_{n}g$,
for some $g\in G$ and some $\phi_{n}$ so,\begin{equation}
\lim_{\lambda\rightarrow0}\phi(\lambda)\cdot A=g^{-1}\lim_{\lambda\rightarrow0}\phi_{\mathsf{n}}(\lambda)\cdot(g\cdot A).\label{eq:limit}\end{equation}
 Writing $g\cdot A$ as\[
g\cdot A=\left(\begin{array}{cc}
a(g) & b(g)\\
c(g) & d(g)\end{array}\right),\]
we obtain \[
\phi_{n}(\lambda)\cdot(g\cdot A)=\left(\begin{array}{ll}
a(g) & b(g)\lambda^{2n}\\
c(g)\lambda^{-2n} & d(g)\end{array}\right).\]
 By the Hilbert-Mumford criterion, the limit (\ref{eq:limit}) exists
for some 1PS, so we must have $c(g)=0$, for some $g\in G$. This
means that $g\cdot A$ is an upper triangular $n$-matrix.
\end{proof}
Note that this result can be easily generalized to describe the stable
points of the action of $SL(m,\mathbb{C})$ under simultaneous conjugation
on the vector space of $n$-tuples of $m\times m$ matrices. To summarize,
for the action of $SL(2,\mathbb{C})$ on $V_{n}$ we have shown the
following, which is a consequence of Theorem \ref{thm:UT} and Proposition
\ref{pro:stable}.

\begin{thm}
\label{thm:0}The following are equivalent for an element $A=(A_{1},...,A_{n})\in V_{n}$,
$n\geq1$.

(i) $A\in V_{n}$ is stable.

(ii) There exists $1\leq j,k,l\leq n$ such that $(A_{j},A_{k},A_{l})\in V_{3}$
is stable.

(iii) There exists $1\leq j,k,l\leq n$ such that $\sigma_{jk}(A)\neq0$
or $\Delta_{jkl}(A)\neq0$.

(iv) $A$ is not similar to an upper triangular $n$-matrix.

(v) There is no proper nonzero subspace of $\mathbb{C}^{2}$ preserved
by the set $\{ A_{1},...,A_{n}\}$.$\hfill\square$
\end{thm}

\section{Irreducibility of representations of finitely generated groups}

We now use the numerical condition for stability found above to derive
a similar criterion for irreducibility of a representation of a finitely
generated group $\Gamma$ into $G=SL(2,\mathbb{C})$. As in the introduction,
by fixing a set of generators $\bar{\epsilon}_{1},...,\bar{\epsilon}_{n}$
of $\Gamma$, we associate to a representation $\rho:\Gamma\to G$
the point $A=(A_{1},...,A_{n})\in X_{n}$ given by\[
A=(\rho(\bar{\epsilon}_{1}),...,\rho(\bar{\epsilon}_{n})).\]

In this section we are therefore dealing with the conjugation action
of $G=SL(2,\mathbb{C})$ restricted to the affine subvariety $X_{n}=G^{\times n}\subset V_{n}$.
For $G$-invariant functions on $X_{n}$, we continue to use the same
notations described in Definition \ref{def:notation}, in particular,
we still denote by $t_{j}$ (resp. $t_{jk}$) trace of the matrix
$A_{j}$ (resp. $A_{j}A_{k}$). Because of the standard identities\begin{eqnarray*}
\mathsf{tr}(B_{1}B_{2})+\mathsf{tr}(B_{1}^{-1}B_{2}) & = & \mathsf{tr}(B_{1})\mathsf{tr}(B_{2})\\
\mathsf{tr}(B_{1}) & = & \mathsf{tr}(B_{1}^{-1})\\
\mathsf{tr}(B_{1}^{2}) & = & \mathsf{tr}^{2}(B_{1})-2,\end{eqnarray*}
valid for any two $SL(2,\mathbb{C})$ matrices $B_{1},B_{2}$, some
of those $G$ invariant functions acquire a new form\begin{eqnarray*}
\nu_{j} & = & \frac{t_{j}^{2}}{2}-2,\\
\sigma_{jk} & = & t_{j}^{2}+t_{k}^{2}+t_{jk}^{2}-t_{j}t_{k}t_{jk}-4=\mathsf{tr}([A_{j},A_{k}])-2,\end{eqnarray*}
where $[B_{1},B_{2}]=B_{1}B_{2}B_{1}^{-1}B_{2}^{-1}$ denotes the
commutator of two $SL(2,\mathbb{C})$ matrices. 

Recall that a representation $\rho:\Gamma\to G$ is called irreducible
if there are no proper nonzero subspaces of $\mathbb{C}^{2}$ which
are invariant under $\rho(\Gamma)$. Therefore, if all matrices $\rho(\gamma)$,
$\gamma\in\Gamma$ are upper triangular, $\rho$ is reducible. In
this case, $\mathsf{tr}([\rho(h_{1}),\rho(h_{2})])=2$ for any $h_{1},h_{2}\in\Gamma$.
This condition is also sufficient for reducibility as proved in \cite{CS}.

\begin{thm}
\emph{(Culler \& Shalen \cite{CS})} \label{thm:CS}A representation
$\rho:\Gamma\to G$ is irreducible if and only if $\mathsf{tr}(c)\neq2$
for some element $c$ of the commutator subgroup $[\Gamma,\Gamma]$.
\end{thm}
Making use of Theorem \ref{thm:0}, we obtain necessary and sufficient
conditions for irreducibility, depending only on the $n$-matrix $A$
associated to $\rho$, via the fixed choice of generators of $\Gamma$.

\begin{thm}
\label{thm: irred}Let $A=(\rho(\bar{\epsilon}_{1}),...,\rho(\bar{\epsilon}_{n}))\in X_{n}$
be the $n$-tuple associated with the representation $\rho:\Gamma\to G$.
Then, the following sentences are equivalent.

(i) $A\in X_{n}$ is stable.

(ii) There exists $1\leq j,k,l\leq n$ such that $(\rho(\bar{\epsilon}_{j}),\rho(\bar{\epsilon}_{k}),\rho(\bar{\epsilon}_{l}))\in X_{3}$
is stable.

(iii) There exists $1\leq j,k,l\leq n$ such that $\sigma_{jk}(A)\neq0$
or $\Delta_{jkl}(A)\neq0$.

(iv) $A$ is not similar to an upper triangular $n$-matrix.

(v) $\rho:\Gamma\rightarrow G$ is irreducible.
\end{thm}
\begin{proof}
The equivalence between (i)-(iv) follows easily from the equivalence
of (i)-(iv) for the $SL(2,\mathbb{C})$ action on $V_{n}$ proved
in Theorem \ref{thm:0}. Let us show that (iv) is equivalent to (v).
If $A\notin\mathcal{UT}_{n}$ there is no subspace of $\mathbb{C}^{2}$
preserved by the set $\{\rho(\bar{\epsilon}_{1}),...,\rho(\bar{\epsilon}_{n})\}$,
so $\rho$ is irreducible. Conversely, if $A\in\mathcal{UT}_{n}$
then, after conjugating $\rho$ with some $g\in G$, all the matrices
$\rho(h)$, $h\in\Gamma$ will be upper triangular because the $\bar{\epsilon}_{j}$
are the generators of $\Gamma$. So $\rho$ is reducible.
\end{proof}
This result completes the proof of theorems \ref{thm:1} and \ref{thm:2};
it can be viewed as a sharpening of Theorem \ref{thm:CS} (it also
generalizes prop. 1.5.5 of \cite{CS}), since part (iii) implies that
irreducibility of a representation $\rho:\Gamma\to G$ can be verified
by computing the values of a finite number (precisely $\binom{n}{2}+\binom{n}{3}=\frac{n^{3}-n}{6}$)
of functions of the $n$ matrices $\rho(\bar{\epsilon}_{j})\in SL(2,\mathbb{C})$,
$j=1,...,n$. Again, note that Theorem \ref{thm: irred} is true for
any algebraically closed field of characteristic 0.

Above, the criterion for irreducibility is written in terms of conditions
for pairs and for triples of $SL(2,\mathbb{C})$ matrices. However,
when working with representations of a free group $F_{n}$, the most
important conditions are the ones for pairs because of the next result,
which also shows that the condition $\nu_{j}=0$ can be easily removed
in the irreducible case.

\begin{prop}
Let $\rho:F_{n}\to G$ be an irreducible representation. Then, there
exists a choice of generators of $F_{n}$ such that the corresponding
$n$-matrix satisfies $\sigma_{12}\neq0$ and $\nu_{1}\neq0$. 
\end{prop}
\begin{proof}
First note that choosing a new set of generators of $F_{n}$ is equivalent
to performing an automorphism of $F_{n}$. Hence, we will find such
an automorphism which, upon acting on the $n$-matrix $A$ associated
with the irreducible $\rho$, will verify $\sigma_{12}\neq0$ and
$\nu_{1}\neq0$. If some $\sigma_{jk}\neq0$ we can just permute the
indices to obtain $\sigma_{12}\neq0$. So, suppose that all $\sigma_{jk}$
are zero and let $(A_{j},A_{k},A_{l})$ be a stable triple, so that
$\Delta_{jkl}\neq0$. Permute the generators again so that the triple
becomes $(A_{1},A_{2},A_{3})$ and assume this triple is already in
the form (\ref{eq:form}) with $e_{1}b_{2}c_{3}\neq0$ and $e_{2}e_{3}+b_{2}c_{3}=0$.
Then, perform the shift automorphism $\epsilon_{1}\mapsto\epsilon_{1}\epsilon_{3}$
of $F_{n}$, which corresponds to $A_{1}\mapsto A_{1}'=A_{1}A_{3}$
and an easy computation shows $\sigma_{12}\mapsto\sigma_{12}':=d_{1}b_{2}c_{3}(d_{1}b_{2}c_{3}+(a_{1}a_{3}-d_{1}d_{3})e_{2})$.
Since none of the factors $d_{1},b_{2},c_{3},e_{2}$ and $a_{3}$
can be zero, the condition $e_{2}e_{3}+b_{2}c_{3}=0$ implies that
$\sigma_{12}'\neq0$. Similarly, using the shift automorphism $A_{1}\mapsto A_{1}A_{2}^{k}$,
for some $k\in\mathbb{Z}$, we end up with $\nu_{1}\neq0$.
\end{proof}
Consider now the stable quotient $X_{n}^{\textrm{st}}/G$, where $X_{n}^{\textrm{st}}\subset X_{n}$
is the subset of stable points in $X_{n}$. As mentioned in the introduction,
this is a geometric quotient and has the structure of an affine algebraic
variety. Because of the identification between $X_{n}^{\textrm{st}}/G$
and $Hom(F_{n},G)^{\textrm{irr}}/G$, where $Hom(F_{n},G)^{\mathrm{irr}}$
is the subset of irreducible representations in $Hom(F_{n},G)$, it
is not difficult to show that $X_{n}^{\textrm{st}}/G$ is nonsingular,
and is therefore a complex analytic manifold of dimension $3n-3$
(see for example \cite{Gu}). Also, when the finitely generated group
$\Gamma=\pi_{1}S$ is the fundamental group of a surface $S$ of genus
$g>1$, the space of conjugacy classes of irreducible representations
$Hom(\pi_{1}S,G)^{\textrm{irr}}/G$ can be given the structure of
a complex manifold (of dimension $6g-6$), and this can be interpreted
as the space of irreducible flat $SL(2,\mathbb{C})$ bundles on $S$
(\cite{Gu}).

In general, when the finitely generated group $\Gamma$ is the fundamental
group of a manifold $M$, it is clear that a representation $\rho:\Gamma\to G$
will define a flat rank 2 vector bundle $E_{\rho}$ over $M$ with
trivial determinant. If $M$ is an algebraic variety, one can consider
moduli spaces for these bundles, and in particular, the moduli spaces
of \emph{stable} and \emph{semistable} vector bundles on compact Riemann
surfaces are well known. 

In this context, one might ask what is the relationship between stability
of the holomorphic vector bundle $E_{\rho}$ and the stability of
the $n$-matrix $A\in X_{n}$ associated with $\rho$. This relation
is simple in one direction. It is a general fact that if $E_{\rho}$
is stable then $\rho$ is irreducible, so $A$ is stable, by what
we saw above. However, the converse is not true and there are irreducible
representations giving rise to unstable vector bundles. An example
of such a vector bundle on a genus $g>1$ Riemann surface $S$ is
the following. Let $L$ be a degree $g-1$ line bundle whose square
is the canonical bundle of $S$. Then the unique (up to isomorphism)
indecomposable vector bundle which is an extension of the form\[
0\to L\to E\to L^{-1}\to0,\]
is associated to an irreducible representation $\rho:\pi_{1}S\to SL(2,\mathbb{C})$.
Actually, $\rho$ is a Schottky representation, as it factors through
a representation $\tilde{\rho}:F_{g}\to G$ of a free group of rank
$g$, for a certain natural projection $\pi_{1}S\to F_{g}$, and $\tilde{\rho}$
defines a Schottky group in $PSL(2,\mathbb{C})$ uniformizing $S$
(see \cite{Fl}). Therefore, the corresponding matrix $A$ is stable,
although $E_{\rho}=E$ is clearly unstable as a vector bundle. Let
$Hom(\pi_{1}S,G)^{\textrm{st}}/G$ be the space of conjugacy classes
of representations $\rho$ such that $E_{\rho}$ is a stable vector
bundle on $S$. Since under the map $\rho\mapsto E_{\rho}$, the space
$Hom(\pi_{1}S,G)^{\textrm{irr}}/G$ is the parameter space for a holomorphic
family of vector bundles over $S$, by a result of Narasimhan and
Seshadri (see \cite{NS}, Thm. 3), we conclude the following.

\begin{prop}
\label{pro:st-irr}The complement of $Hom(\pi_{1}S,G)^{\textrm{st}}/G$
inside $Hom(\pi_{1}S,G)^{\textrm{irr}}/G$ is a nonempty analytic
subset. 
\end{prop}
The characterization of this analytic subset in terms of the geometry
of $S$ seems to be a difficult open problem.

\section{The Magnus trace map $T_{n}$}

In this section, we will study the Magnus trace map $T_{n}$ defined
in the introduction\begin{eqnarray*}
T_{n}:X_{n} & \rightarrow & \mathbb{C}^{3n-3}\\
A & \mapsto & (t_{1},t_{2},t_{12},t_{3},t_{13},t_{23},...,t_{k},t_{1k},t_{2k},...,t_{n},t_{1n},t_{2n})\end{eqnarray*}
and determine when a given point in $\mathbb{C}^{3n-3}$ determines
a finite (and nonzero) number of orbits of the action of $G=SL(2,\mathbb{C})$
on $X_{n}=G^{\times n}$. Recall the definitions (\ref{eq:FD}) of
$\sigma_{jk}$ and of the Fricke discriminant $\Delta_{jkl}$. One
of the results in \cite{M} states that

\begin{thm}
\emph{(Magnus \cite{M})} Let $n\geq4$ and $A_{1},A_{2},A_{3}\in SL(2,\mathbb{C})$
be three fixed matrices verifying $\sigma_{12}\neq0$ and $\Delta_{123}\neq0$.
Let $\mathbf{p}=(t_{1},t_{2},t_{12},t_{3},t_{13},t_{23})\in\mathbb{C}^{6}$.
Then given any $\mathbf{q}\in\mathbb{C}^{3n-9}$ there exist $A_{4},...,A_{n}$
such that \[
T_{n}(A_{1},...,A_{n})=(\mathbf{p},\mathbf{q}).\]
The number of solutions $(A_{4},...,A_{n})$ is bounded by $2^{n-2}$.
$\hfill\square$
\end{thm}
This result could suggest that the image of $T_{n}$ does not intersect
the sets in the image corresponding to the conditions $\sigma_{12}=0$
and $\Delta_{123}=0$, as these are defined by polynomials in the
variables $t_{j}$ and $t_{jk}$. It turns out that the Fricke discriminant
condition is not really necessary, and we only need to prevent $\sigma_{12}$
from being zero. 

\begin{thm}
\label{thm:surject}Let $n\geq2$ and $A_{1},A_{2}\in SL(2,\mathbb{C})$
be fixed matrices with $\sigma_{12}\ne0$, and let $\mathbf{r}=(t_{1},t_{2},t_{12})$.
Then given any $\mathbf{s}\in\mathbb{C}^{3n-6}$ there exist $A_{3},...,A_{n}$
such that \[
T_{n}(A_{1},...,A_{n})=(\mathbf{r},\mathbf{s}).\]

\end{thm}
In proving Theorem \ref{thm:surject}, we will use systematically
the transposition invariant forms of a pair of matrices $A_{1},A_{2}$
verifying $\sigma_{12}\neq0$, given in Proposition \ref{pro:ST}.
We will also regard elements of $SL(2,\mathbb{C})$ as complexified
$SU(2)$ matrices, or as complexified quaternions with unit norm,
writing them in the form\begin{equation}
A_{j}=\left(\begin{array}{cc}
\alpha_{j}+i\beta_{j} & \gamma_{j}+i\delta_{j}\\
-\gamma_{j}+i\delta_{j} & \alpha_{j}-i\beta_{j}\end{array}\right)\in SL(2,\mathbb{C}),\quad j=1,...,n,\label{eq:SU2C}\end{equation}
with\begin{equation}
\alpha_{j},\beta_{j},\gamma_{j},\delta_{j}\in\mathbb{C}\textrm{, and}\quad\alpha_{j}^{2}+\beta_{j}^{2}+\gamma_{j}^{2}+\delta_{j}^{2}=1.\label{eq:norma}\end{equation}
For matrices in $G=SL(2,\mathbb{C})$, and using the parameterizations
of $A_{j}$ as in (\ref{eq:SU2C}), the following is an immediate
consequence of the proof of Proposition \ref{pro:ST}.

\begin{prop}
\label{pro:ST2}Let $A=(A_{1},A_{2})\in X_{2}$ verify $\sigma_{12}\neq0$.
If $\nu_{1}\neq0$, then $A$ is similar to a pair $B=(B_{1},B_{2})$
of the form\[
B_{1}=\left(\begin{array}{cc}
\alpha_{1}+i\beta_{1} & 0\\
0 & \alpha_{1}-i\beta_{1}\end{array}\right),\qquad B_{2}=\left(\begin{array}{cc}
\alpha_{2}+i\beta_{2} & i\delta_{2}\\
i\delta_{2} & \alpha_{2}-i\beta_{2}\end{array}\right),\]
with $\alpha_{1}^{2}+\beta_{1}^{2}=1$, $\alpha_{2}^{2}+\beta_{2}^{2}+\delta_{2}^{2}=1$
and $\delta_{2}^{2}=\frac{\sigma_{12}}{2\nu_{1}}\neq0$. In case $\nu_{1}=0$
and $\nu_{2}\neq0$, the situation is completely analogous switching
$B_{1}$ with $B_{2}$. If $\nu_{1}=\nu_{2}=0$, then $A$ is similar
to a pair $B=(B_{1},B_{2})$ of the form\[
B_{1}=\left(\begin{array}{cc}
\alpha_{1}+\lambda & i\lambda\\
i\lambda & \alpha_{1}-\lambda\end{array}\right),\qquad B_{2}=\left(\begin{array}{cc}
\alpha_{2}-\lambda & i\lambda\\
i\lambda & \alpha_{2}+\lambda\end{array}\right),\]
with $\alpha_{1}^{2}=\alpha_{2}^{2}=1$ and $\lambda^{4}=\frac{\sigma_{12}}{16}\neq0$.
$\hfill\square$
\end{prop}
To adapt our notation to the trace map $T_{n}$, let us denote the
components of an arbitrary element $\mathbf{z}$ of the range $\mathbb{C}^{3n-3}$
by\begin{equation}
\mathbf{z}=(z_{1},z_{2},z_{12},...,z_{k},z_{1k},z_{2k},...,z_{n},z_{1n},z_{2n}),\label{eq:z}\end{equation}
and define\begin{eqnarray*}
\nu_{j}(\mathbf{z}) & = & \frac{z_{j}^{2}}{2}-2,\\
\sigma_{jk}(\mathbf{z}) & = & z_{j}^{2}+z_{k}^{2}+z_{jk}^{2}-z_{j}z_{k}z_{jk}-4.\end{eqnarray*}
In this way, given any polynomial $p(\mathbf{z})$ in the coordinates
of $\mathbf{z}$, its pullback under $T_{n}$ will be the same polynomial
in the variables $(t_{1},t_{2},t_{12},...,t_{n},t_{1n},t_{2n})$.
To simplify the notation, we will sometimes write these polynomial
functions without reference to the variables. Consider the following
Zariski closed subsets of $X_{n}$ and of $\mathbb{C}^{3n-3}$, respectively,
\begin{eqnarray*}
U_{12} & = & \{ A\in X_{n}:\sigma_{12}(A)=0\}\\
Z_{12} & = & \{\mathbf{z}\in\mathbb{C}^{3n-3}:\sigma_{12}(\mathbf{z})=0\},\end{eqnarray*}
so that $T_{n}(U_{12})=Z_{12}$. It is immediate that Theorem \ref{thm:surject}
above can be restated as

\begin{thm}
\label{thm:surject2}The algebraic map $T_{n}:X_{n}\setminus U_{12}\rightarrow\mathbb{C}^{3n-3}\setminus Z_{12}$
is surjective for $n\geq2$.
\end{thm}
\begin{proof}
We need to show that, for any given $\mathbf{z}$ as in (\ref{eq:z})
verifying $z_{1}^{2}+z_{2}^{2}+z_{12}^{2}-z_{1}z_{2}z_{12}\neq4$,
there is a set of $n$ matrices \[
A_{j}=\left(\begin{array}{cc}
\alpha_{j}+i\beta_{j} & \gamma_{j}+i\delta_{j}\\
-\gamma_{j}+i\delta_{j} & \alpha_{j}-i\beta_{j}\end{array}\right)\in SL(2,\mathbb{C}),\quad j=1,...,n,\]
 such that $\mathsf{tr}(A_{j})=z_{j}$ and $\mathsf{tr}(A_{j}A_{k})=z_{jk}$,
for all coordinates of $\mathbf{z}$. Let us start with the case when
$\nu_{1}(\mathbf{z})\neq0$ which means $z_{1}\neq\pm2$ and corresponds
under $T_{n}$ to $\nu_{1}(A)\neq0$. Then, using Proposition \ref{pro:ST2},
if there is a solution, there is one in the form\begin{equation}
A_{1}=\left(\begin{array}{cc}
\alpha_{1}+i\beta_{1} & 0\\
0 & \alpha_{1}-i\beta_{1}\end{array}\right),\qquad A_{2}=\left(\begin{array}{cc}
\alpha_{2}+i\beta_{2} & i\delta_{2}\\
i\delta_{2} & \alpha_{2}-i\beta_{2}\end{array}\right),\label{eq:geral1}\end{equation}
\[
A_{k}=\left(\begin{array}{cc}
\alpha_{k}+i\beta_{k} & \gamma_{k}+i\delta_{k}\\
-\gamma_{k}+i\delta_{k} & \alpha_{k}-i\beta_{k}\end{array}\right),\quad k=3,...,n.\]
Then, we need to solve\begin{eqnarray}
z_{k} & = & 2\alpha_{k},\qquad k=1,...,n,\nonumber \\
z_{1k} & = & \mathsf{tr}(A_{1}A_{k})=2(\alpha_{1}\alpha_{k}-\beta_{1}\beta_{k}),\qquad k=2,...,n,\label{eq:Tngeral}\\
z_{2k} & = & \mathsf{tr}(A_{2}A_{k})=2(\alpha_{2}\alpha_{k}-\beta_{2}\beta_{k}-\delta_{2}\delta_{k}),\qquad k=3,...,n,\nonumber \end{eqnarray}
and an explicit solution is obtained by setting $\alpha_{k}=z_{k}/2$,
for $k=1,...,n$ and \[
\beta_{1}=\sqrt{1-\alpha_{1}^{2}},\qquad\beta_{k}=\frac{z_{1}z_{k}-\frac{1}{2}z_{1k}}{\beta_{1}},\qquad k=2,...,n,\]
\begin{equation}
\delta_{2}=\sqrt{q_{22}},\quad\delta_{k}=\frac{q_{2k}}{\delta_{2}},\qquad k=3,...,n,\quad\textrm{and}\label{eq:geral2}\end{equation}
\[
\gamma_{k}=\frac{\sqrt{q_{22}q_{kk}-q_{2k}^{2}}}{\delta_{2}},\qquad k=3,...,n,\]
 where we have used the abbreviations $q_{kk}=1-\alpha_{k}^{2}-\beta_{k}^{2}$
and $q_{jk}=\alpha_{j}\alpha_{k}-\beta_{j}\beta_{k}-\frac{1}{2}z_{jk}$,
if $j\neq k$. The denominators $\beta_{1}$ and $\delta_{2}$ are
both nonzero because of our assumptions on $\nu_{1}(\mathbf{z})=-2\beta_{1}^{2}$
and $\sigma_{12}(\mathbf{z})=-4\beta_{1}^{2}\delta_{2}^{2}=2\nu_{1}\delta_{2}^{2}$.
Observe that, after assuming $A_{1}$ diagonal and $A_{2}=A_{2}^{T}$,
the above are the only solutions, and different solutions correspond
to different choices of square roots in the above expressions. Also
note that the Fricke discriminant appears as $\Delta_{12k}=16\beta_{1}^{2}\delta_{2}^{2}\gamma_{k}^{2}$. 

When $\nu_{1}(\mathbf{z})=\nu_{2}(\mathbf{z})=0$, from Proposition
\ref{pro:ST2}, we may use the parametrization\begin{equation}
A_{1}=\left(\begin{array}{cc}
\alpha_{1}+\lambda & i\lambda\\
i\lambda & \alpha_{1}-\lambda\end{array}\right),\qquad A_{2}=\left(\begin{array}{cc}
\alpha_{2}-\lambda & i\lambda\\
i\lambda & \alpha_{2}+\lambda\end{array}\right),\qquad\lambda\neq0\label{eq:espec1}\end{equation}
\[
A_{k}=\left(\begin{array}{cc}
\alpha_{k}+i\beta_{k} & \gamma_{k}+i\delta_{k}\\
-\gamma_{k}+i\delta_{k} & \alpha_{k}-i\beta_{k}\end{array}\right),\quad k=3,...,n,\]
 where $\alpha_{1}$ and $\alpha_{2}$ are square roots of 1. Now,
the equations $T_{n}(A)=\mathbf{z}$ become\begin{eqnarray}
z_{k} & = & 2\alpha_{k}\nonumber \\
z_{1k} & = & 2\alpha_{1}\alpha_{k}+2\lambda(i\beta_{k}-\delta_{k})\label{eq:espec2}\\
z_{2k} & = & 2\alpha_{2}\alpha_{k}+2\lambda(-i\beta_{k}-\delta_{k}),\qquad k=3,...,n\nonumber \end{eqnarray}
 which one can easily solve for $\alpha_{k}$, $\beta_{k}$, and $\delta_{k}$.
The variables $\gamma_{k}$ are then obtained from the normalization
(\ref{eq:norma}). As before, the solutions in this form will be in
finite number.
\end{proof}
As remarked, the proof of Theorem \ref{thm:surject2} implies that
any point $\mathbf{z}\in\mathbb{C}^{3n-3}\setminus Z_{12}$ determines
an orbit up to a finite ambiguity. More concretely, we describe all
the solutions as follows.

\begin{thm}
Let $\mathbf{z}\in\mathbb{C}^{3n-3}\setminus Z_{12}$ and $A=(A_{1},A_{2},...,A_{n})$
be a solution of $T_{n}(A)=\mathbf{z}$, such that $A_{1}$ and $A_{2}$
are invariant under transposition. Then $T_{n}^{-1}(z)$ is the $G$
orbit of the finite set\begin{equation}
\left\{ (A_{1},A_{2},B_{3},...,B_{n}):\  B_{k}=A_{k}\textrm{ or }B_{k}=A_{k}^{T}\right\} \subset X_{n}\label{eq:orbit}\end{equation}
of cardinality $\leq2^{n-2}$.
\end{thm}
\begin{proof}
First, let $\mathbf{z}\in\mathbb{C}^{3n-3}\setminus Z_{12}$ verify
$z_{1}\neq\pm2$, and let $\mathcal{S}_{1}$ be the space of solutions
of $T_{n}(A)=\mathbf{z}$ with $A_{1}$ diagonal and $A_{2}^{T}=A_{2}$.
From the proof of Theorem (\ref{thm:surject2}) and from Propoposition
(\ref{pro:ST2}), any $A\in\mathcal{S}_{1}$ is obtained from formulas
(\ref{eq:geral1}) and (\ref{eq:geral2}), for a certain choice of
square roots of $\beta_{1}$, $\delta_{2}$ and $\gamma_{3}$,...,
$\gamma_{n}$. Changing all the signs of $\beta_{1},...,\beta_{n},\gamma_{3},...,\gamma_{n}$
simultaneously, gives a well defined map $\sigma_{1}:\mathcal{S}_{1}\to\mathcal{S}_{1}$.
Similarly, let $\sigma_{2}:\mathcal{S}_{1}\to\mathcal{S}_{1}$ be
the operation of changing simultaneously the signs of $\delta_{2},...,\delta_{n},\gamma_{3},...,\gamma_{n}$.
The signs of each $\gamma_{k}$, $k=3,...,n$ can be changed independently
of the rest, giving maps $\sigma_{k}:\mathcal{S}_{1}\to\mathcal{S}_{1}$,
$k=3,...,n$. For matrices in the form (\ref{eq:SU2C}), it is clear
that conjugation by $g=(0,i,i,0)$ changes the triple of vectors $(\beta,\gamma,\delta)$
into $(-\beta,-\gamma,\delta)$, and conjugation by $g=(i,0,0,-i)$
maps $(\beta,\gamma,\delta)$ into $(\beta,-\gamma,-\delta)$ (naturally
any conjugation fixes the vector $\alpha$). Therefore, $\sigma_{1}$
and $\sigma_{2}$ act trivially on the space $\mathcal{S}_{1}/G$.
Since transposition of a single $A_{k}$ changes the sign of a single
$\gamma_{k}$, $k=3,...,n$, the orbits are as in (\ref{eq:orbit}).
Similarly, to treat the case $z_{1}=\pm2$, let $\mathcal{S}_{0}$
be the space of solutions of $T_{n}(A)=\mathbf{z}$ of the form (\ref{eq:espec1}),
(\ref{eq:espec2}). As before, define $\sigma_{0}:\mathcal{S}_{0}\to\mathcal{S}_{0}$
by changing simultaneously the signs of $\lambda,\beta_{1},...,\beta_{n},\delta_{1},...,\delta_{n}$.
By the same reason as before, this acts trivially on $\mathcal{S}_{0}/G$,
so also in this case, the set of solutions is given by (\ref{eq:orbit}).
\end{proof}
We have thus finished the proof of Theorem \ref{thm:Magnus}. 

\appendix

\section{Reconstruction of matrices from traces}

In this Appendix, we briefly indicate the generic reconstruction of
(the $SL(2,\mathbb{C})$ orbit of) general $n$-tuples of $2\times2$
matrices $A\in V_{n}$ from a minimal number of traces given by the
analogous trace map\begin{eqnarray}
\hat{T}_{n}:V_{n} & \rightarrow & \mathbb{C}^{4n-3}\nonumber \\
A & \mapsto & (t_{1},t_{11},t_{2},t_{22},t_{12},t_{3},t_{33},t_{13},t_{23},...,t_{n},t_{nn},t_{1n},t_{2n}).\label{eq:TjVn}\end{eqnarray}
The formulas are identical to the ones in Theorem \ref{thm:surject2},
when written only in terms of the functions $\tau_{jk}:\mathbb{C}^{4n-3}\to\mathbb{C}$,
given by \[
\tau_{jk}(\mathbf{z})=z_{jk}-\frac{1}{2}z_{j}z_{k},\qquad j,k\in\{1,n\},\]
 except that now we allow coordinates $z_{jk}$ with $j=k$ and, writing
$A_{j}$ in the form (\ref{eq:SU2C}) we don't require the normalization
condition (\ref{eq:norma}). For instance, for $\sigma_{12}(\mathbf{z})\neq0$,
and $\nu_{1}(\mathbf{z})\neq0$, the following matrices form a solution
to $\hat{T}_{n}(A)=\mathbf{z}$.\[
A_{1}=\left(\begin{array}{cc}
\alpha_{1}+i\beta_{1} & 0\\
0 & \alpha_{1}-i\beta_{1}\end{array}\right),\qquad A_{2}=\left(\begin{array}{cc}
\alpha_{2}+i\beta_{2} & i\delta_{2}\\
i\delta_{2} & \alpha_{2}-i\beta_{2}\end{array}\right),\]
\[
A_{k}=\left(\begin{array}{cc}
\alpha_{k}+i\beta_{k} & \gamma_{k}+i\delta_{k}\\
-\gamma_{k}+i\delta_{k} & \alpha_{k}-i\beta_{k}\end{array}\right),\quad k=3,...,n,\]
with\[
\beta_{1}=\sqrt{-\frac{\tau_{11}}{2}}\qquad\beta_{k}=-\frac{\tau_{1k}}{2\beta_{1}},\qquad k=2,...,n,\]
\[
\delta_{2}=\sqrt{\frac{\sigma_{12}}{2\tau_{11}}},\qquad\delta_{k}=\frac{\tau_{2k}\tau_{11}-\tau_{12}\tau_{1k}}{2\delta_{2}\tau_{11}},\qquad,k=3,...,n,\]
\[
\gamma_{k}=\sqrt{-\frac{\Delta_{12k}}{4\sigma_{12}}},\qquad k=3,...,n,\]
where the `Fricke discriminant', as a function of $\tau_{jk}(\mathbf{z})$,
is given by\[
\Delta_{12k}(\mathbf{z})=2(\tau_{12}^{2}\tau_{kk}+\tau_{1k}^{2}t_{22}+\tau_{2k}^{2}\tau_{11}-2\tau_{12}\tau_{1k}\tau_{2k}-\tau_{11}\tau_{22}\tau_{kk}),\]
in agreement with (\ref{eq:det}). The case with $z_{1},z_{2}\in\{-2,2\}$
can be treated similarly.

\section{Examples with $T_{n}^{-1}(\mathbf{z})=\emptyset$.}

Here, we show that the image of $T_{n}$, for $n\geq4$, does not
contain certain points $\mathbf{z}\in\mathbb{C}^{3n-3}$. From Theorem
\ref{thm:surject}, if the equation $T_{n}(A)=\mathbf{z}$ has no
solution $A\in X_{n}$, then $\sigma_{12}(\mathbf{z})=0$. Take, for
example, any $\mathbf{z}$ verifying $z_{1}=z_{2}=z_{12}=2$, and
for which there are two indices $k,l\notin\{1,2\}$ such that\begin{equation}
(z_{1k}-z_{k})(z_{2l}-z_{l})\neq(z_{1l}-z_{l})(z_{2k}-z_{k}).\label{eq:hypot}\end{equation}
 Then, without loss of generality $A$ can be taken in the form:\[
A_{1}=\left(\begin{array}{cc}
1 & b_{1}\\
0 & 1\end{array}\right),\qquad A_{2}=\left(\begin{array}{cc}
1 & b_{2}\\
0 & 1\end{array}\right),\]
\[
A_{j}=\left(\begin{array}{cc}
\alpha_{j}+i\beta_{j} & b_{j}\\
c_{j} & \alpha_{j}-i\beta_{j}\end{array}\right),\quad j=3,...,n.\]
Suppose $T_{n}(A)=\mathbf{z}$. Then, we should have, for all $j=3,...,n$,\begin{eqnarray}
z_{j} & = & 2\alpha_{j}\nonumber \\
z_{1j} & = & 2\alpha_{j}+b_{1}c_{j}\label{eq:counter}\\
z_{2j} & = & 2\alpha_{j}+b_{2}c_{j},\nonumber \end{eqnarray}
In the case that $b_{1}$ and $b_{2}$ are nonzero, this implies \[
c_{k}=\frac{z_{1k}-z_{k}}{b_{1}}=\frac{z_{2k}-z_{k}}{b_{2}},\qquad c_{l}=\frac{z_{1l}-z_{l}}{b_{1}}=\frac{z_{2l}-z_{l}}{b_{2}},\]
which is impossible by our hypothesis (\ref{eq:hypot}). Similarly,
if $b_{1}=0$, equations (\ref{eq:counter}) imply $z_{1k}=z_{k}$
and $z_{1l}=z_{l}$, contradicting again (\ref{eq:hypot}). The same
happens if $b_{2}=0$. 

More generally, a similar situation occurs in the case that $\sigma_{12}(\mathbf{z})=0$
but, for example, $\nu_{1}(\mathbf{z})$ is nonzero. In this case,
without loss of generality (see lemma \ref{lemma}), we have the parametrization\[
A_{1}=\left(\begin{array}{cc}
\alpha_{1}+i\beta_{1} & 0\\
0 & \alpha_{1}-i\beta_{1}\end{array}\right),\]

\[
A_{2}=\left(\begin{array}{cc}
\alpha_{2}+i\beta_{2} & 1\\
0 & \alpha_{2}-i\beta_{2}\end{array}\right),\quad\textrm{or}\quad A_{2}=\left(\begin{array}{cc}
\alpha_{2}+i\beta_{2} & 0\\
0 & \alpha_{2}-i\beta_{2}\end{array}\right),\]

\[
A_{k}=\left(\begin{array}{cc}
\alpha_{k}+i\beta_{k} & b_{k}\\
c_{k} & \alpha_{k}-i\beta_{k}\end{array}\right),\quad k=3,...,n,\]
 Assuming $A_{2}$ non-diagonal, we need to solve the equations, for
every $k=3,...,n$. \begin{eqnarray}
z_{k} & = & 2\alpha_{k}\nonumber \\
z_{1k} & = & \mathsf{tr}(A_{1}A_{k})=2(\alpha_{1}\alpha_{k}-\beta_{1}\beta_{k})\nonumber \\
z_{2k} & = & \mathsf{tr}(A_{2}A_{k})=2(\alpha_{2}\alpha_{k}-\beta_{2}\beta_{k})+c_{k}.\label{eq:desapare}\end{eqnarray}
 One finds that $\beta_{1}=\sqrt{-\frac{\nu_{1}(\mathbf{z})}{2}}\neq0$,
$\beta_{2}=\sqrt{-\frac{\nu_{2}(\mathbf{z})}{2}}$ and\begin{eqnarray*}
\beta_{k} & = & \frac{2\alpha_{1}\alpha_{k}-z_{1k}}{2\beta_{1}}=-\frac{\tau_{1k}(\mathbf{z})}{2\beta_{1}}\\
c_{k} & = & \frac{\beta_{1}\tau_{2k}(\mathbf{z})-\beta_{2}\tau_{1k}(\mathbf{z})}{\beta_{1}}\end{eqnarray*}
 If $c_{k}\neq0$, then the $b_{k}$ can be found using the normalization
$\alpha_{k}^{2}+\beta_{k}^{2}-b_{k}c_{k}=1.$ However, when $\beta_{1}\tau_{2k}(\mathbf{z})=\beta_{2}\tau_{1k}(\mathbf{z})$,
for some $k\neq1,2$, in the above expression we obtain $c_{k}=0$
which implies that $A_{2}$ has to be in diagonal form, and then $c_{k}$
disappears from equation (\ref{eq:desapare}), so we need to have
$\beta_{1}\tau_{2j}(\mathbf{z})=\beta_{2}\tau_{1j}(\mathbf{z})$ for
all $j\neq1,2$. Obviously, this does not hold in general. Note also
that these computations, combined with the previous case $z_{1}=\pm2$,
allows one to reprove surjectivity for $n=3$ (in this case, $c_{3}=0$
is exactly equivalent to the necessary condition $\beta_{1}\tau_{23}(\mathbf{z})=\beta_{2}\tau_{13}(\mathbf{z})$).

\end{document}